# On the output-input stability property for multivariable nonlinear control systems

Daniel Liberzon

Coordinated Science Laboratory
University of Illinois at Urbana-Champaign
liberzon@uiuc.edu

November 16, 2018

## Abstract

We study the recently introduced notion of output-input stability, which is a robust variant of the minimum-phase property for general smooth nonlinear control systems. The subject of this paper is developing the theory of output-input stability in the multi-input, multi-output setting. We show that output-input stability can be viewed as a combination of two system properties, one related to detectability and the other to left-invertibility. For systems affine in controls, we provide a necessary and sufficient condition for output-input stability, which relies on Hirschorn's nonlinear structure algorithm.

## 1 Introduction

For systems with inputs, two properties of interest are asymptotic stability under zero inputs and bounded state response to bounded inputs. It is well known that for linear time-invariant systems the first property implies the second one, but for nonlinear systems this is not the case. The notion of *input-to-state stability* (ISS) introduced in [11] captures both of the above properties; it requires that bounded inputs produce bounded states and inputs converging (or equal) to zero produce states converging to zero.

Dual concepts of detectability result if one considers systems with outputs. For linear systems, one of equivalent ways to define detectability is to demand that the state converge to zero along every trajectory for which the output is identically zero. The notion of *output-to-state stability* (OSS) introduced in [12] is a robust version of the detectability property for nonlinear systems and a dual of ISS; it requires that the state be bounded if the output is bounded and converge to zero if the output converges to zero.

The present line of work is concerned with the *minimum-phase* property of systems with both inputs and outputs. A linear system is minimum-phase if whenever the output is identically zero, both the state and the input must converge to zero; in the frequency domain, this is characterized by stability of system zeros. Byrnes and Isidori [1] provided an important and natural extension of the minimum-phase property to nonlinear systems (affine in controls). According to their definition, a system is minimum-phase if its *zero dynamics*—the internal dynamics of the system under the action of an input that holds the output constantly at zero—are asymptotically stable.

The above remarks suggest that to complete the picture, one should have a robust version of the last property, which should ask the state and the input to be bounded when the output is bounded and to become small when the output is small. Such a concept was proposed in the recent paper [5] under the name of *output-input stability*. In the spirit of ISS, output-input stability requires the state and the input of the system to be bounded by a suitable function of the output and derivatives of the output, modulo a decaying term depending on initial conditions. This property is in general stronger than the minimum-phase property[1] defined in [1]. Output-input stability can be studied with the help of the tools that have been developed over the years to study ISS, OSS, and related notions (such as Lyapunov-like dissipation inequalities); the minimum-phase property, on the other hand, is investigated using different techniques (such as computation of zero dynamics and normal forms). This makes output-input stability an appropriate alternative notion to use in those situations where the minimum-phase property is insufficient or difficult to check.

The results of [5] provide a fairly complete theory of output-input stable single-input, single-output (SISO) nonlinear control systems. In this paper we continue to study the output-input stability property for multi-input, multi-output (MIMO) systems. Our goal is to investigate a connection between output-input stability and structural properties of control systems which have been studied in the context of system inversion. In particular, we show the relevance of Hirschorn's nonlinear structure algorithm [2] (which is an extension of Silverman's linear structure algorithm [8, 9]) in establishing output-input stability. After providing necessary definitions and preliminary results in Section 2 and reviewing the nonlinear structure algorithm in Section 3, we derive our main result for affine systems in Section 4 and then discuss some extensions in Section 5.

## 2 Background and preliminary results

Consider the system

$$\begin{aligned}\dot{x} &= f(x,u) \\ y &= h(x)\end{aligned} \qquad (1)$$

where the state $x$ takes values in $\mathbb{R}^n$, the input $u$ takes values in $\mathbb{R}^m$, the output $y$ takes values in $\mathbb{R}^p$ (for some positive integers $n$, $m$, and $p$), and the functions $f$ and $h$ are smooth. In this paper we restrict admissible input (or "control") signals to be at least continuous. For every initial condition $x(0)$ and every input $u(\cdot)$, there is a solution $x(\cdot)$ of (1) defined on a maximal interval $[0, T_{\max})$, and the corresponding output $y(\cdot)$. We write $\mathcal{C}^k$ for the space of $k$ times continuously differentiable functions $u : [0, \infty) \to \mathbb{R}^m$, where $k$ is some nonnegative integer. Whenever the input $u$ is in $\mathcal{C}^k$, the derivatives $\dot{y}, \ddot{y}, \ldots, y^{(k+1)}$ exist and are continuous; they are given by

$$y^{(i)}(t) = H_i\bigl(x(t), u(t), \ldots, u^{(i-1)}(t)\bigr), \qquad i = 1, \ldots, k+1, \quad t \in [0, T_{\max}) \qquad (2)$$

where for $i = 0, 1, \ldots$ the functions $H_i : \mathbb{R}^n \times (\mathbb{R}^m)^i \to \mathbb{R}^p$ are defined recursively via $H_0 := h$ and

$$H_{i+1}(x, u_0, \ldots, u_i) := \frac{\partial H_i}{\partial x} f(x, u_0) + \sum_{j=0}^{i-1} \frac{\partial H_i}{\partial u_j} u_{j+1}$$

---

[1]Strictly speaking, this statement only makes sense for systems affine in controls, because otherwise the minimum-phase property is not defined. For example, the scalar system $\dot{y} = 1 + y^2 + u^2$ is output-input stable (because $|u| \leq \sqrt{\dot{y}}$) but not minimum-phase (in fact, no input can hold the output at zero).

(here the arguments of $H_i$ are $x \in \mathbb{R}^n$ and $u_0, \ldots, u_{i-1} \in \mathbb{R}^m$). Given a nonnegative integer $k$, we will denote by $\mathbf{y}^k$ the $\mathbb{R}^{p(k+1)}$-valued signal

$$\mathbf{y}^k := \begin{pmatrix} y \\ \dot{y} \\ \vdots \\ y^{(k)} \end{pmatrix}$$

provided that the indicated derivatives exist. We will let $\|\cdot\|_{[a,b]}$ denote the supremum norm of a signal restricted to an interval $[a,b]$, i.e., $\|z\|_{[a,b]} := \sup\{|z(s)| : a \leq s \leq b\}$, where $|\cdot|$ is the standard Euclidean norm.

According to Definition 1 of [5], the system (1) is called *output-input stable* if there exist a positive integer $N$, a class $\mathcal{KL}$ function[2] $\beta$, and a class $\mathcal{K}_\infty$ function $\gamma$ such that for every initial state $x(0)$ and every input $u \in \mathcal{C}^{N-1}$ the inequality

$$\left|\begin{pmatrix} x(t) \\ u(t) \end{pmatrix}\right| \leq \beta(|x(0)|, t) + \gamma(\|\mathbf{y}^N\|_{[0,t]}) \tag{3}$$

holds for all $t$ in the domain of the corresponding solution. (The assumption that $u$ belongs to $\mathcal{C}^{N-1}$ is made to guarantee that $\mathbf{y}^N$ is well defined, and can be weakened if the function $H_N$ is independent of $u_{N-1}$.) As discussed in [5], the concept of output-input stability finds applications in feedback stabilization, adaptive control, and other areas.

It is perhaps best to interpret output-input stability as a combination of two separate properties of the system. The first one is expressed by the inequality

$$|x(t)| \leq \beta(|x(0)|, t) + \gamma(\|\mathbf{y}^N\|_{[0,t]}) \tag{4}$$

and corresponds to detectability (OSS) with respect to the output and its derivatives, uniform over inputs; this property was called "weak uniform 0-detectability of order $N$" in [5]. The results of [4, 12] imply that the system (1) is weakly uniformly 0-detectable of order $N$ if there exists a continuously differentiable, positive definite, radially unbounded function $V : \mathbb{R}^n \to \mathbb{R}$ that satisfies

$$\frac{\partial V}{\partial x} f(x, u_0) \leq -\alpha(|x|) + \chi\left(\left|\begin{pmatrix} H_0(x) \\ \ldots \\ H_N(x, u_0, \ldots, u_{N-1}) \end{pmatrix}\right|\right) \quad \forall\, x, u_0, \ldots, u_{N-1}$$

for some functions $\alpha, \chi \in \mathcal{K}_\infty$. As discussed in [5], the class of weakly uniformly 0-detectable systems includes all affine systems in global normal form with ISS internal dynamics.

The second ingredient of the output-input stability property is described by the inequality

$$|u(t)| \leq \beta(|x(0)|, t) + \gamma(\|\mathbf{y}^N\|_{[0,t]}) \tag{5}$$

which says that the input should become small if the output and its derivatives are small. Loosely speaking, this suggests that the system has a stable left inverse in the input-output sense. Unlike uniform detectability, this property does not seem to admit a Lyapunov-like characterization. In the SISO case it is closely related to the existence of a relative degree; see [5, Theorem 1]. In general, however, this second property needs to be understood better, which is precisely the goal of the present paper. In what follows, we formulate and study a useful property which, in combination with (4), yields (5).

---

[2] Recall that a function $\alpha : [0, \infty) \to [0, \infty)$ is said to be of *class* $\mathcal{K}$ if it is continuous, strictly increasing, and $\alpha(0) = 0$. If $\alpha \in \mathcal{K}$ is unbounded, then it is said to be of *class* $\mathcal{K}_\infty$. A function $\beta : [0, \infty) \times [0, \infty) \to [0, \infty)$ is said to be of *class* $\mathcal{KL}$ if $\beta(\cdot, t)$ is of class $\mathcal{K}$ for each fixed $t \geq 0$ and $\beta(s, t)$ decreases to 0 as $t \to \infty$ for each fixed $s \geq 0$.

Let us say that the system (1) has the *input-bounding property* if there exist a positive integer $k^*$ and two class $\mathcal{K}_\infty$ functions $\rho_1$ and $\rho_2$ such that we have

$$|u_0| \leq \rho_1(|x|) + \rho_2 \left( \left| \begin{pmatrix} H_0(x) \\ \ldots \\ H_{k^*}(x, u_0, \ldots, u_{k^*-1}) \end{pmatrix} \right| \right) \qquad \forall\, x, u_0, \ldots, u_{k^*-1}. \tag{6}$$

(We could omit the first component $H_0(x)$ of the argument of $\rho_2$ without loss of generality, but will keep it for notational convenience.) Defined in this way, the input-bounding property represents a functional relation between the input and state variables. The next result recasts this property in terms of trajectories of the system.

**Lemma 1** *The system* (1) *has the input-bounding property if and only if there exist a positive integer $k^*$ and two class $\mathcal{K}_\infty$ functions $\rho_1$ and $\rho_2$ such that for every initial condition and every input $u \in \mathcal{C}^{k^*-1}$ the inequality*

$$|u(t)| \leq \rho_1(|x(t)|) + \rho_2(|\mathbf{y}^{k^*}(t)|) \tag{7}$$

*holds for all t in the domain of the corresponding solution.*

PROOF. In view of (2), it is clear that (6) implies (7), with the same $k^*$, $\rho_1$, and $\rho_2$. To show the converse, suppose that (6) is violated for some $x, u_0, \ldots, u_{k^*-1}$. Take $x$ to be the initial condition and apply an input $u$ satisfying $u^{(i)}(0) = u_i$, $i = 0, \ldots, k^* - 1$. Then it is easy to see that (7) does not hold for small $t$. □

We point out that the input-bounding property resembles in its appearance the notion of relative degree as defined in [5], but is actually much less restrictive (especially for MIMO systems). The next result reveals the connection between output-input stability, weak uniform 0-detectability, and the input-bounding property; it was essentially proved in [5].

**Lemma 2** *Suppose that the system* (1) *is weakly uniformly 0-detectable of order $\hat{k}$:*

$$|x(t)| \leq \bar{\beta}(|x(0)|, t) + \bar{\gamma}(\|\mathbf{y}^{\hat{k}}\|_{[0,t]}) \tag{8}$$

*where $\hat{k} \geq 0$, $\bar{\beta} \in \mathcal{KL}$, and $\bar{\gamma} \in \mathcal{K}_\infty$. Suppose also that it has the input-bounding property as defined above. Then it is output-input stable, with $N = \max\{\hat{k}, k^*\}$.*

PROOF. Since (1) has the input-bounding property, by Lemma 1 the inequality (7) holds with $\rho_1, \rho_2 \in \mathcal{K}_\infty$. Combining this with (8) and using the simple fact that for every class $\mathcal{K}$ function $\rho$ and arbitrary numbers $s_1, s_2 \geq 0$ one has $\rho(s_1 + s_2) \leq \rho(2s_1) + \rho(2s_2)$, we arrive at the inequality (3) with $N := \max\{\hat{k}, k^*\}$, $\beta(s,t) := \rho_1(2\bar{\beta}(s,t)) + \bar{\beta}(s,t)$, and $\gamma(s) := \rho_1(2\bar{\gamma}(s)) + \rho_2(s) + \bar{\gamma}(s)$. Thus (1) is output-input stable. □

The above observation explains the importance of the input-bounding property. As we will see later, a natural way of checking this property for systems affine in controls is provided by Hirschorn's nonlinear structure algorithm.

# 3  Nonlinear structure algorithm

In this section and the next one we assume that $m \leq p$ and the system (1) takes the form

$$\dot{x} = f(x) + G(x)u \qquad (9)$$
$$y = h(x)$$

Its dynamics can also be written in more detail as

$$\dot{x} = f(x) + \sum_{i=1}^{m} g_i(x) u_i.$$

We assume that $f(0) = 0$ and $h(0) = 0$ (although the second assumption is only made for convenience and can be removed). For $f$ a vector field and $R$ a scalar-valued function, $L_f R$ denotes the derivative of $R$ along $f$. We will also let $L_f$ act on vectors and matrices componentwise; for example, if $R$ is matrix-valued, then $L_f R$ is a matrix-valued function whose $ij$th component is given by $L_f(R_{ij})$. All functions are assumed to have the smoothness required for all relevant derivatives to exist. Dimensions of vectors and matrices will be omitted when clear from the context.

To make the paper self-contained, we now describe (a global version of) Hirschorn's structure algorithm from [2]. It iteratively constructs a sequence of vectors $z_k \in \mathbb{R}^p$ which take the form

$$z_k = M_k(x) \mathbf{y}^k$$

(i.e., their components are $x$-dependent linear combinations of the original outputs and their derivatives) and satisfy

$$z_k = h_k(x) + J_k(x) u$$

where

$$J_k(x) = \begin{pmatrix} \bar{J}_k(x) \\ 0 \end{pmatrix}$$

with $\bar{J}_k(x)$ an $r_k \times m$ matrix of constant rank $r_k \geq 0$. Partitioning all vectors accordingly, we can rewrite this as

$$\begin{pmatrix} \bar{z}_k \\ \hat{z}_k \end{pmatrix} = \begin{pmatrix} \overline{M}_k(x) \\ \widehat{M}_k(x) \end{pmatrix} \mathbf{y}^k = \begin{pmatrix} \bar{h}_k(x) \\ \hat{h}_k(x) \end{pmatrix} + \begin{pmatrix} \bar{J}_k(x) \\ 0 \end{pmatrix} u. \qquad (10)$$

The initial data for $k = 0$ are given by $M_0 \equiv I$, $\bar{z}_0 = 0$, $\hat{z}_0 = y$, $\bar{h}_0 \equiv 0$, $\hat{h}_0 = h$, $\bar{J}_0 \equiv 0$, and $r_0 = 0$. To obtain $z_{k+1}$, we differentiate $\hat{z}_k$ along solutions of the system. On one hand, we have

$$\dot{\hat{z}}_k = \frac{d}{dt}(\widehat{M}_k(x) \mathbf{y}^k) = L_f \widehat{M}_k(x) \mathbf{y}^k + \sum_{i=1}^{m} L_{g_i} \widehat{M}_k(x) u_i \mathbf{y}^k + \widehat{M}_k(x) \dot{\mathbf{y}}^k.$$

ASSUMPTION 1. $L_{g_i} \widehat{M}_k(x) \equiv 0$ for all $i$.

If this assumption holds, then the dependence on $u$ disappears and we obtain

$$\dot{\hat{z}}_k = L_f \widehat{M}_k(x) \mathbf{y}^k + \widehat{M}_k(x) \dot{\mathbf{y}}^k. \qquad (11)$$

On the other hand,

$$\dot{z}_k = \frac{d}{dt}\hat{h}_k(x) = L_f\hat{h}_k(x) + \sum_{i=1}^{m} L_{g_i}\hat{h}_k(x)u_i = L_f\hat{h}_k(x) + \widehat{J}_k(x)u$$

where

$$\widehat{J}_k(x) := \frac{\partial \hat{h}_k}{\partial x}G(x).$$

We thus arrive at the equation

$$\begin{pmatrix} \bar{z}_k \\ \dot{\hat{z}}_k \end{pmatrix} = \begin{pmatrix} \bar{h}_k(x) \\ L_f\hat{h}_k(x) \end{pmatrix} + \begin{pmatrix} \bar{J}_k(x) \\ \widehat{J}_k(x) \end{pmatrix} u. \tag{12}$$

ASSUMPTION 2. *The matrix* $\begin{pmatrix} \bar{J}_k(x) \\ \widehat{J}_k(x) \end{pmatrix}$ *has a constant rank* $r_{k+1}$ *and a fixed set of* $r_{k+1}$ *rows that are linearly independent for all* $x$.

This assumption implies that there exists a $p \times p$ permutation matrix $E_k$ (independent of $x$) such that

$$E_k \begin{pmatrix} \bar{J}_k(x) \\ \widehat{J}_k(x) \end{pmatrix} = \begin{pmatrix} \bar{J}_{k+1}(x) \\ \widetilde{J}_{k+1}(x) \end{pmatrix}$$

where $\bar{J}_{k+1}(x)$ is an $r_{k+1} \times m$ matrix of rank $r_{k+1}$. (In fact, typically it suffices to reorder the rows of $\widehat{J}_k(x)$, since the rows of $\bar{J}_k(x)$ are already linearly independent.) Therefore, a $(p - r_{k+1}) \times r_{k+1}$ matrix $F_k(x)$ can be found such that $F_k(x)\bar{J}_{k+1}(x) + \widetilde{J}_{k+1}(x) \equiv 0$. Consider the $p \times p$ matrix

$$R_k(x) := \begin{pmatrix} I & 0 \\ F_k(x) & I \end{pmatrix} E_k.$$

Multiplying both sides of (12) by $R_k(x)$, we obtain

$$R_k(x) \begin{pmatrix} \bar{z}_k \\ \dot{\hat{z}}_k \end{pmatrix} = R_k(x) \begin{pmatrix} \bar{h}_k(x) \\ L_f\hat{h}_k(x) \end{pmatrix} + \begin{pmatrix} \bar{J}_{k+1}(x) \\ 0 \end{pmatrix} u.$$

The iteration step is completed by setting

$$z_{k+1} := R_k(x) \begin{pmatrix} \bar{z}_k \\ \dot{\hat{z}}_k \end{pmatrix}, \qquad h_{k+1}(x) := R_k(x) \begin{pmatrix} \bar{h}_k(x) \\ L_f\hat{h}_k(x) \end{pmatrix}$$

and noticing that from (10) and (11) we have $z_{k+1} = M_{k+1}(x)\mathbf{y}^{k+1}$ where

$$M_{k+1}(x) = \begin{pmatrix} \overline{M}_k(x) & 0 \\ \left(L_f\widehat{M}_k(x) \quad 0\right) + \left(0 \quad \widehat{M}_k(x)\right) \end{pmatrix}.$$

It is easy to see that $r_{k+1} \geq r_k$ for all $k$. If for some $k$ we have $r_k = m$, then the algorithm terminates. We denote the smallest such $k$ by $k^*$ or let $k^* = \infty$ if the above condition is never met.

It is not hard to show that Assumption 1 is fulfilled if $L_{g_i}L_f^j F_k(x) \equiv 0$ for $0 \leq j + k \leq k^* - 2$ and $1 \leq i \leq m$ [2, Theorem 2]. (In particular, it is trivially satisfied if $k^* = 1$.) Hirschorn also gives sufficient conditions for this in terms of the original system [2, Corollary 1]. Assumption 2 represents a global regularity condition. We discuss possible ways of relaxing these assumptions in Section 5. The following example illustrates the application of the above algorithm to a system satisfying both assumptions at each step.

EXAMPLE 1. Consider the system

$$\begin{aligned}
\dot{x}_1 &= u_1 \\
\dot{x}_2 &= x_3 + x_4 u_1 \\
\dot{x}_3 &= u_2 \\
\dot{x}_4 &= -x_4 + x_1^2 \\
y &= (x_1, x_2)^T
\end{aligned} \tag{13}$$

We have

$$\begin{pmatrix} \dot{y}_1 \\ \dot{y}_2 \end{pmatrix} = \begin{pmatrix} 0 \\ x_3 \end{pmatrix} + \begin{pmatrix} 1 & 0 \\ x_4 & 0 \end{pmatrix} \begin{pmatrix} u_1 \\ u_2 \end{pmatrix}$$

so $r_1 = 1$. Multiplying both sides by $R_0 = \begin{pmatrix} 1 & 0 \\ -x_4 & 1 \end{pmatrix}$ yields

$$\begin{pmatrix} \dot{y}_1 \\ -x_4 \dot{y}_1 + \dot{y}_2 \end{pmatrix} = \begin{pmatrix} 0 \\ x_3 \end{pmatrix} + \begin{pmatrix} 1 & 0 \\ 0 & 0 \end{pmatrix} \begin{pmatrix} u_1 \\ u_2 \end{pmatrix}.$$

Differentiating the bottom component, we obtain

$$\begin{pmatrix} \dot{y}_1 \\ (x_4 - x_1^2)\dot{y}_1 - x_4 \ddot{y}_1 + \ddot{y}_2 \end{pmatrix} = \begin{pmatrix} 1 & 0 \\ 0 & 1 \end{pmatrix} \begin{pmatrix} u_1 \\ u_2 \end{pmatrix} \tag{14}$$

hence $r_2 = 2$ and the algorithm terminates with $k^* = 2$. $\square$

## 4 Main result

We are now ready to state and prove our main result for affine systems.

**Theorem 1** *Let Assumptions 1 and 2 hold for each $k \geq 0$. Then the system (9) is output-input stable if and only if it is weakly uniformly 0-detectable (of some order $\hat{k}$) and the algorithm of Section 3 gives $k^* < \infty$.*

PROOF. We first prove sufficiency. If $k^* < \infty$, then the system is left-invertible, i.e., one can explicitly solve for $u$ as a function of $x$ and $\mathbf{y}^{k^*}$. Indeed, we have

$$\overline{M}_{k^*}(x)\mathbf{y}^{k^*} = \bar{h}_{k^*}(x) + \bar{J}_{k^*}(x)u$$

where $\bar{J}_{k^*}(x)$ is square and invertible. Multiplying both sides by $\bar{J}_{k^*}^{-1}(x)$, we obtain

$$u = -\bar{J}_{k^*}^{-1}(x)\bar{h}_{k^*}(x) + \bar{J}_{k^*}^{-1}(x)\overline{M}_{k^*}(x)\mathbf{y}^{k^*} =: A(x) + B(x)\mathbf{y}^{k^*}.$$

It is easy to show that $f(0) = 0$ implies $A(0) = 0$. It follows that

$$|A(x)| \leq \gamma_1(|x|), \qquad |B(x)| \leq |B(0)| + \gamma_2(|x|)$$

for some $\gamma_1, \gamma_2 \in \mathcal{K}_\infty$. This yields

$$|u| \leq \gamma_1(|x|) + \big(|B(0)| + \gamma_2(|x|)\big)|\mathbf{y}^{k^*}| \leq \gamma_1(|x|) + |B(0)||\mathbf{y}^{k^*}| + \frac{1}{2}(\gamma_2(|x|))^2 + \frac{1}{2}|\mathbf{y}^{k^*}|^2$$

and so the inequality (7) holds with $\rho_1(r) := \gamma_1(r) + \frac{1}{2}(\gamma_2(r))^2$ and $\rho_2(r) := |B(0)|r + \frac{1}{2}r^2$. By Lemma 1 the system has the input-bounding property, and Lemma 2 guarantees output-input stability with $N = \max\{\hat{k}, k^*\}$.

Necessity follows from the correspondence between Hirschorn's structure algorithm and the zero dynamics algorithm, which is explained, e.g., in [7, Sections 11.1–11.2] or [3, Section 6.1]. If $k^* = \infty$, then in some neighborhood of the origin there exist a vector $\alpha(x) \in \mathbb{R}^m$ and an $m \times l$ matrix $\beta(x)$, where $0 < l \leq m$, such that every feedback law of the form $u = \alpha(x) + \beta(x)v$, $v \in \mathbb{R}^l$ holds the output constantly at zero (for a proper choice of initial conditions). Clearly, this violates the inequality (5), hence the system cannot be output-input stable. The fact that output-input stability implies weak uniform 0-detectability is an immediate consequence of the definitions. □

EXAMPLE 1 (continued). Consider again the system (13). From the formula (14) we conclude that $|u_1| \leq |\dot{y}_1|$ and
$$|u_2| \leq \frac{1}{2}(x_4 - x_1^2)^2 + \frac{1}{2}\dot{y}_1^2 + \frac{1}{2}x_4^2 + \frac{1}{2}\ddot{y}_1^2 + |\ddot{y}_2|$$
hence the system has the input-bounding property. It is also weakly uniformly 0-detectable of order 1, as is seen from the bound $|x_3| = |\dot{y}_2 - x_4\dot{y}_1| \leq |\dot{y}_2| + |x_4||\dot{y}_1|$ and the fact that the equation for $x_4$ is ISS with respect to $x_1$. Therefore, the system (13) is output-input stable. □

A sufficient condition for termination of Hirschorn's algorithm, and consequently for the input-bounding property, is the existence of a uniform (vector) relative degree as defined, e.g., in [3]. This condition is by no means necessary; note that the system considered in the above example does not have a uniform relative degree.

The next example illustrates what can happen when Assumption 2 is violated.

EXAMPLE 2. Consider the system

$$\begin{aligned}
\dot{x}_1 &= u_1 \\
\dot{x}_2 &= x_3 + x_2 u_2 \\
\dot{x}_3 &= u_2 \\
\dot{x}_4 &= -x_4 + x_1^2 \\
y &= (x_1, x_2)^T
\end{aligned} \qquad (15)$$

We have
$$\begin{pmatrix} \dot{y}_1 \\ \dot{y}_2 \end{pmatrix} = \begin{pmatrix} 0 \\ x_3 \end{pmatrix} + \begin{pmatrix} 1 & 0 \\ 0 & x_2 \end{pmatrix} \begin{pmatrix} u_1 \\ u_2 \end{pmatrix}$$
and the rank of the matrix on the right-hand side drops from 2 to 1 when $x_2 = 0$. It is not difficult to show that we can pick a bounded sequence of initial states along which $x_2(0)$ converges to 0, a sequence of values of $u_2(0)$ converging to $\infty$, and appropriately chosen sequences of values for $\dot{u}_2(0)$, $\ddot{u}_2(0)$, … such that the derivatives $\dot{y}_2(0)$, $\ddot{y}_2(0)$, … are all kept at zero. Also, let $u_1 \equiv 0$ so that $y_1 \equiv 0$. This implies that the inequality (5) is violated for small $t$, hence the system (15) is not output-input stable. (The proof of Theorem 1 in [5] contains a general argument along these lines.) □

It is instructive to note that both the system (13) considered in Example 1 and the system (15) considered in Example 2 are minimum-phase, with zero dynamics in both cases being given by $\dot{x}_4 = -x_4$. An important fact not elucidated by the zero dynamics is that the minimum-phase property of the system (13) is "robust" (small $y$, $\dot{y}$, … force $x$ and $u$ to be small) while the minimum-phase property of the system (15) is "fragile" (small $y$, $\dot{y}$, … can correspond to arbitrarily large $u$).

# 5 Some extensions

Assumptions 1 and 2 are automatically satisfied for SISO affine systems with uniform relative degree in the sense of [3]. In the SISO case, output-input stability actually implies the existence of relative degree for a class of systems which includes systems affine in controls; see [5, Theorem 1]. We will see below that in contrast with that result, for MIMO systems Assumptions 1 and 2 are not necessary for output-input stability.

Suppose that at some step $\bar{k}$ Assumption 1 does not hold. Then we can use the following modification of Hirschorn's algorithm, proposed in [10]. Let us collect all $u$-dependent terms on the right-hand side:

$$L_f \widehat{M_{\bar{k}}}(x)\mathbf{y}^{\bar{k}} + \widehat{M_{\bar{k}}}(x)\dot{\mathbf{y}}^{\bar{k}} = L_f \hat{h}_{\bar{k}}(x) + \sum_{i=1}^m L_{g_i} \hat{h}_{\bar{k}}(x) u_i - \sum_{i=1}^m L_{g_i} \widehat{M_{\bar{k}}}(x) u_i \mathbf{y}^{\bar{k}} =: L_f \hat{h}_{\bar{k}}(x) + \widehat{J_{\bar{k}}}(x, \mathbf{y}^{\bar{k}}) u.$$

Then we can follow the same steps as before and define matrices $R_k$, $k \geq \bar{k}$ and $J_k$, $k > \bar{k}$ which will depend rationally on $\mathbf{y}^k$, and vectors $z_k$, $k > \bar{k}$ which will depend rationally but not necessarily linearly on $\mathbf{y}^k$. Assuming that Assumption 2 is still satisfied with "for all $x, \mathbf{y}^{k}$" replacing "for all $x$", and that $k^* < \infty$, we can bound $u$ in terms of $x$ and a rational function of $\mathbf{y}^{k^*}$ (cf. [7, Example 11.20]). This in general does not give the input-bounding property. However, in many cases the above dependance on $\mathbf{y}^{k^*}$ is actually polynomial, and then the input-bounding property does hold. The next example serves to support this observation.

EXAMPLE 3. Consider the system

$$\begin{aligned}
\dot{x}_1 &= u_1 \\
\dot{x}_2 &= x_3 + x_2 u_1 \\
\dot{x}_3 &= u_2 \\
y &= (x_1, x_2)^T
\end{aligned} \tag{16}$$

We have

$$\begin{pmatrix} \dot{y}_1 \\ \dot{y}_2 \end{pmatrix} = \begin{pmatrix} 0 \\ x_3 \end{pmatrix} + \begin{pmatrix} 1 & 0 \\ x_2 & 0 \end{pmatrix} \begin{pmatrix} u_1 \\ u_2 \end{pmatrix}$$

so $r_1 = 1$. Multiplying both sides by $R_0 = \begin{pmatrix} 1 & 0 \\ -x_2 & 1 \end{pmatrix}$ gives

$$\begin{pmatrix} \dot{y}_1 \\ -x_2 \dot{y}_1 + \dot{y}_2 \end{pmatrix} = \begin{pmatrix} 0 \\ x_3 \end{pmatrix} + \begin{pmatrix} 1 & 0 \\ 0 & 0 \end{pmatrix} \begin{pmatrix} u_1 \\ u_2 \end{pmatrix}.$$

Differentiating the bottom component, we obtain

$$\begin{pmatrix} \dot{y}_1 \\ (-x_3 - x_2 u_1)\dot{y}_1 - x_2 \ddot{y}_1 + \ddot{y}_2 \end{pmatrix} = \begin{pmatrix} 1 & 0 \\ 0 & 1 \end{pmatrix} \begin{pmatrix} u_1 \\ u_2 \end{pmatrix}.$$

Moving the $u_1$-dependent term from the left-hand side to the right, we can rewrite this as

$$\begin{pmatrix} \dot{y}_1 \\ -x_3 \dot{y}_1 - x_2 \ddot{y}_1 + \ddot{y}_2 \end{pmatrix} = \begin{pmatrix} 1 & 0 \\ x_2 \dot{y}_1 & 1 \end{pmatrix} \begin{pmatrix} u_1 \\ u_2 \end{pmatrix}.$$

Finally, multiplying both sides by $\begin{pmatrix} 1 & 0 \\ x_2 \dot{y}_1 & 1 \end{pmatrix}^{-1} = \begin{pmatrix} 1 & 0 \\ -x_2 \dot{y}_1 & 1 \end{pmatrix}$, we arrive at

$$\begin{pmatrix} \dot{y}_1 \\ -x_2 \dot{y}_1^2 - x_3 \dot{y}_1 - x_2 \ddot{y}_1 + \ddot{y}_2 \end{pmatrix} = \begin{pmatrix} u_1 \\ u_2 \end{pmatrix}$$

and now the input-bounding property can be shown as in Example 1. Weak uniform 0-detectability can also be verified along the same lines as in Example 1. Therefore, the system (16) is output-input stable. We leave it to the reader to check that if we add another state to the system by replacing the equation $\dot{x}_3 = u_2$ with $\dot{x}_3 = x_4$ and adding the equation $\dot{x}_4 = u_2$, then the conclusion remains unchanged. □

Note that Assumption 1 was not used in the proof of the necessity part of Theorem 1. It is not hard to see that necessity is preserved if this assumption is removed and the structure algorithm is modified as described above.

We now demonstrate that a system can be output-input stable even if Assumption 2 is not satisfied. (This is to be contrasted with Example 2.)

EXAMPLE 4. The system

$$\dot{x}_1 = u_1$$
$$\dot{x}_2 = x_5 + x_4 u_2$$
$$\dot{x}_3 = x_4$$
$$\dot{x}_4 = u_2$$
$$\dot{x}_5 = u_3$$
$$y = (x_1, x_2, x_3)^T$$

is output-input stable, as can be seen from the formulas $u_1 = \dot{y}_1$, $u_2 = \ddot{y}_3$, $u_3 = \ddot{y}_2 - \ddot{y}_3^2 - \dot{y}_3 \dddot{y}_3$, $x_4 = \dot{y}_3$, and $x_5 = \dot{y}_2 - \dot{y}_3 \ddot{y}_3$. However, when we try to apply Hirschorn's algorithm, we obtain

$$\begin{pmatrix} \dot{y}_1 \\ \dot{y}_2 \\ \dot{y}_3 \end{pmatrix} = \begin{pmatrix} 0 \\ x_5 \\ x_4 \end{pmatrix} + \begin{pmatrix} 1 & 0 & 0 \\ 0 & x_4 & 0 \\ 0 & 0 & 0 \end{pmatrix} \begin{pmatrix} u_1 \\ u_2 \\ u_3 \end{pmatrix}$$

and the matrix on the right-hand side does not have a constant rank[3]. □

It is straightforward to obtain a local variant of Theorem 1. Namely, if Assumption 2 holds for all $k$ only in a neighborhood of a given equilibrium—say, the origin—and if in the definition of output-input stability we only require the inequality (3) to be satisfied for initial conditions close to the origin and for small $t$, then the theorem is still valid. (Note also that locally, the second statement in Assumption 2 automatically follows from the first.) For analytic systems, Assumption 2 is guaranteed to hold on an open and dense subset of $\mathbb{R}^n$ (the complement of a singular set). Local versions of output-input stability may be useful for control design purposes away from singularities (cf. [13]).

The above results can be used to establish output-input stability of some nonaffine systems. Note that to have the input-bounding property, we only need to be able to bound—and not necessarily solve for—the input in terms of the state and derivatives of the output. As a simple generalization, consider a system of the form

$$\dot{x} = f(x) + \sum_{i=1}^{m} g_i(x) \gamma_i(u_i)$$

where the functions $\gamma_i$, $i = 1, \ldots, m$ are bounded from below by some class $\mathcal{K}_\infty$ functions. It is easy to show that if the associated "virtual input" system

$$\dot{x} = f(x) + \sum_{i=1}^{m} g_i(x) v_i$$

---
[3]This system has trivial zero dynamics ($x \equiv 0$) but $x = 0$ is not a regular point of the zero dynamics algorithm.

is covered by the sufficiency part of Theorem 1, then the original nonaffine system is output-input stable. Of course, left-invertibility of the virtual input system is not necessary. One can even have more inputs than outputs; for example, the scalar system $\dot y = u_1^2 + u_2^4$ is clearly output-input stable.

There has been a lot of activity in the 1980s on establishing invertibility conditions for nonlinear systems in terms of differential geometric concepts—particularly, controlled invariant submanifolds. These and related developments are documented in [3, Chapter 6], [7, Chapter 11], and the references therein, and build on earlier results from geometric theory of linear systems [6]. This body of work helps one to better understand the above constructions and reveals their coordinate independence (see, e.g., the proof of the necessity part of Theorem 1). However, these results seem to be of limited utility in establishing output-input stability, as they do not capture properties such as the one expressed by Assumption 1 or the polynomial dependence discussed earlier in this section. It would be interesting to gain more insight into the questions studied here with the help of geometric tools.